\documentclass{article}

\usepackage[english]{babel}

\usepackage[letterpaper,top=2cm,bottom=2cm,left=3cm,right=3cm,marginparwidth=1.75cm]{geometry}

\usepackage{amsfonts}
\usepackage{graphicx}
\usepackage[all]{xy}
\usepackage{amsmath,textcomp,amssymb,geometry,graphicx,enumerate}
\usepackage{mathrsfs}
\usepackage[colorlinks=true, allcolors=blue]{hyperref}
\title{Counting Non-abelian Coverings of Algebraic Curve}
\author{Peisheng Yu}

\date{}

\begin{document}
\maketitle

\begin{abstract}
In this article, we study the etale coverings of an algebraic curve $C$ with Galois group a semi-direct product $\mathbb{Z}/m\mathbb{Z} \rtimes \mathbb{Z}/n\mathbb{Z}$. Especially, for a given etale cyclic $n$-covering $D \to C$, we determine how many curves $E$ are there, satisfying $E \to D$ is an etale cyclic $m$-covering and $E \to C$ is Galois with non-abelian Galois group, under the assumption $gcd(m,n)=1$.
\end{abstract}

\section{Introduction}
\numberwithin{equation}{section}
Cyclic covering between algebraic varieties is frequently encountered in algebraic geometry. In the case of algebraic curve, the unramified cyclic coverings are encoded by the torsion points of Jacobian. The goal of this short article is to have a glance at the more general coverings between algebraic curves: we study the unramified Galois covering with Galois group equals a semi-direct product $\mathbb{Z}/m\mathbb{Z} \rtimes \mathbb{Z}/n\mathbb{Z}$. Troughout this article, the curves we consider are smooth and projective, defined over complex number $\mathbb{C}$; $m,n \geq 2$ are positive integers. For genus $g=0$, $\mathbb{P}^{1}$ does not have unramified covering other that $\mathbb{P}^{1}$ itself; For genus $g=1$, the fundamental group is isomorphic to $\mathbb{Z} \oplus \mathbb{Z}$, which is abelian. Thus, genus one curve does not have non-abelian unramified covering. For these reasons, the curve $C$ appeared in this article is assumed having genus $g \geq 2$, smooth and projective defined over $\mathbb{C}$. We describe how to construct such coverings $E \to C$ in section 4. Especially, in the case of $gcd(m,n)=1$, we give an explicit formula counting the number of such coverings in section 4 and section 5.

\section{Background: Unramified Cyclic Covering}

Cyclic $n$-coverings of an algebraic variety is well-understood. Roughly speaking, they correspond to the $n$-th roots of a certain line bundle over that variety. For the sake of completeness, we record the well-known facts (\cite{barth2015compact}) concerning the unramified case as in the following subsections.

\subsection{Construct Etale Coverings}

Let $X$ be a smooth projective variety, and assume that there is a free action of $\mathbb{Z}/n \mathbb{Z}$ on $X$. Then the quotient $Y=X/ (\mathbb{Z}/n \mathbb{Z})$ is also smooth and projective. Thus, $\pi: X \to Y$ is an etale cyclic covering of degree $n$.
\\ \hspace{\fill} \\
The pushforward $\pi_{*}\mathcal{O}_{X}$ has a structure of $\mathcal{O}_{Y}$ algebra. Under the action of $\mathbb{Z}/n \mathbb{Z}$, $\pi_{*}\mathcal{O}_{X}$ decomposes into irreducible representations of $\mathbb{Z}/n \mathbb{Z}$. Namely, we have:
\begin{equation}
    \pi_{*}\mathcal{O}_{X} \cong \mathcal{O}_{Y} \oplus \mathcal{L}^{-1} \oplus \cdots \oplus \mathcal{L}^{-n+1} \nonumber
\end{equation}
where $L \in Pic(Y)$ is a $n$-torsion line bundle, which means $\mathcal{L}^{\otimes n}\cong \mathcal{O}_{Y}$.
\\ \hspace{\fill} \\
Conversely, take a $n$-torsion line bundle $\mathcal{L}$ over $Y$, we can construct such an $X$ explicitly: fixing an isomorphism $\mathcal{L}^{\otimes n}\cong \mathcal{O}_{Y}$, we can endow the sheaf $\bigoplus_{i=0}^{n-1}\mathcal{L}^{-i}$ with the structure of a $\mathcal{O}_{Y}$-algebra. Let $X=Spec(\bigoplus_{i=0}^{n-1}\mathcal{L}^{-i})$, then the $\mathcal{O}_{Y}$-structure map $\mathcal{O}_{Y} \to \bigoplus_{i=0}^{n-1}\mathcal{L}^{-i}$ gives a surjection $X \to Y$, which is an etale cyclic $n$-covering. Note that under this construction, in order to guarantee $X$ to be connected, we require $n$ to be minimum: i.e. $n=ord(\mathcal{L})$ in the group $Pic(Y)$. Moreover, if we replace $\mathcal{L}$ by any $\mathcal{L}^{k}$ with $gcd(k,n)=1$, the vairety $X'$ we obtain is isomorphic to $X$ constructed above, because $\bigoplus_{i=0}^{n-1}\mathcal{L}^{-i} \cong \bigoplus_{i=0}^{n-1}\mathcal{L}^{-ki}$ as $\mathcal{O}_{Y}$-algebra. Therefore, if one wants to count $X$ up to isomorphism, such equvilance $\mathcal{L} \sim \mathcal{L}^{k}$ needs to be modulo out.
\\ \hspace{\fill} \\
More geometrically, suppose we have a $n$-torsion line bundle $\mathcal{L}$ on $Y$, and the constant section $1\in \Gamma(Y, \mathcal{O}_{Y})$. We denote $\mathbb{L}$ the total space of $\mathcal{L}$, and we let $p: \mathbb{L} \to Y$ be the bundle projection. If $t\in \Gamma(\mathbb{L}, p^{*}\mathcal{L})$ is the tautological section, the zero divisor of $p^{*}1 - t^{n}$ defines $X$ in $\mathbb{L}$. $\pi = p|_{X}$ exhibits $X$ as an etale cyclic $n$-covering of $Y$. For such $X$, we can also verify that $\pi_{*}\mathcal{O}_{X} \cong \bigoplus_{i=0}^{n-1}\mathcal{L}^{-i} $. Under such description, let $K=\left\langle \sigma \right\rangle \cong \mathbb{Z} / n\mathbb{Z}$, and let $\chi$ denote the generating character of $K$, i.e. $\chi(\sigma)=\zeta_{n}$, a primitive $n$-th root of unity. Then $\mathcal{L}^{-i}$ is the eigenspace of $\chi^{i}$.

\subsection{Case of Curve}

The above argument exhibits a bijection between $n$-torsion points in $Pic(Y)$ and etale cyclic $n$-covering of $Y$. In the case of a smooth projective curve $C$ of genus $g \geq 2$, $n$-torsion points in $Pic^{0}(C) \subset Pic(C)$ is just $J_{C}[n]$, where $J_{C}$ is the Jacobian of $C$. Over $\mathbb{C}$, we know that $J_{C}[n]\cong (\mathbb{Z} / n\mathbb{Z})^{\oplus 2g}$. If we denote $S \subset (\mathbb{Z} / n\mathbb{Z})^{\oplus 2g}$ as the primitive elements in $(\mathbb{Z} / n\mathbb{Z})^{\oplus 2g}$: i.e. $S= \{ v \in (\mathbb{Z} / n\mathbb{Z})^{\oplus 2g}, ord(v)=n \}$, then up to isomorphisms, the connected etale cyclic $n$-covering of $C$ is bijective to $S/ (\mathbb{Z} / n\mathbb{Z})^{\times} $.

\section{Motivation}
This study is motivated by a result considering dihedral covering of $\mathbb{P}^{1}$. Roughly speaking, this result shows that all the etale cyclic $n$-coverings of a hyperelliptical curve give rise to dihedral coverings of $\mathbb{P}^{1}$, as is indicated in the following theorem:
\\ \hspace{\fill} \\
\textbf{Theorem 3.1:} Suppose $H$ is a hyperelliptical curve of genus $g \geq 2$, such that the 2-to-1 covering $H \to \mathbb{P}^{1}$ ramified at $2g+2$ points. There is a hyperelliptical involution $\tau : H \to H$ interchanging the two sheets of $H$. Consider $\pi : X \to H$, which is an etale cyclic $n$-covering of $H$. Then the involution $\tau$ lifts to an involution $\Tilde{\tau}: X \to X$. In other words, $X \to \mathbb{P}^{1}$ is Galois, with Galois group $D_{2n}$.
\\ \hspace{\fill} \\
\textbf{Proof:} Suppose $X$ corresponds to a $n$-torsion line bundle $\mathcal{L}$ over $H$, which by section 2.2 also corresponds to a $n$-torsion point in $J_{H}[n]\cong (\mathbb{Z} / n\mathbb{Z})^{\oplus 2g}$. We need to determine the induced action of $\tau$ on $J_{H}$. Assume the affine piece of $H$ over $A^{1}$ is given by equation $y^{2} =f(x)$, where $f(x)$ is a degree $2g+2$ polynomial in $x$. Then $\omega_{i}=\frac{x^{i-1}dx}{y}, i=1,2,\cdots , g$ form a basis of $H^{0}(H, \omega_{H})$. Since the hyperelliptical involution $\tau$ sends $(x,y)$ to $(x,-y)$, it sends $\omega_{i}$ to $-\omega_{i}$. Since $J_{H} = H^{0}(H, \omega_{H})^{*}/H_{1}(H, \mathbb{Z}) $, we see that $\tau$ induces $[-1]: J_{H} \to J_{H}$. Combining with $J_{H}[n]\cong (\mathbb{Z} / n\mathbb{Z})^{\oplus 2g}$, the induced action of $\tau$ sends $v \in (\mathbb{Z} / n\mathbb{Z})^{\oplus 2g} $ to $-v$. Thus, the torsion line bundle $\mathcal{L}$ is sent to $\mathcal{L}^{-1}$. Note that $X=Spec(\bigoplus_{i=0}^{n-1}\mathcal{L}^{-i})$, we see $\Tilde{\tau}$ sends $X$ to $Spec(\bigoplus_{i=0}^{n-1}\mathcal{L}^{i})\cong X$. Therefore, $\Tilde{\tau}$ is an isomorphism of $X$, satisfying $\Tilde{\tau}^{2}=id_{X}$. From this, we know $X \to \mathbb{P}^{1}$ is Galois.
\\ To determine the Galois group of this etale covering, let $\sigma$ denotes the generator of the Galois group of the covering $X\to H$, which is isomorphic to $\mathbb{Z}/ n \mathbb{Z}$. Recall that $\sigma$ acts on $\bigoplus_{i=0}^{n-1}\mathcal{L}^{-i}$, thus acts on $X$, in the following way: Each $\mathcal{L}^{k}$ is an eigenspace of $\sigma$. Suppose the eigenvalue of $\sigma$ acting on $\mathcal{L}$ is $\zeta$, for some primitive $n$-th root of unity $\zeta$, then the eigenvalue of $\sigma$ acting on $\mathcal{L}^{k}$ equals to $\zeta^{k}$. Together with the fact that $\Tilde{\tau}$ sends $\mathcal{L}^{i}$ to $\mathcal{L}^{-i}$, one easily checks that $\Tilde{\tau} \sigma \Tilde{\tau} = \sigma^{-1}$. Therefore, the Galois group of $X \to \mathbb{P}^{1}$ is $\left\langle \sigma, \Tilde{\tau} | \Tilde{\tau} \sigma \Tilde{\tau} = \sigma^{-1} \right\rangle \cong D_{2n}$. 
\\ \hspace{\fill} \\
Non-abelian covering appears in the theorem above: the dihedral covering. Inspiring by that, we want to construct non-abelian covering for more general curves. Instead of $H \to \mathbb{P}^{1}$, our background setting will be $D \to C$, which is an etale cyclic $n$-covering between smooth projective curves over $\mathbb{C}$, and the genus of $C$ is $g \geq 2$. The goal is to find all curves $E$ up to isomorphism, such that $E \to D$ is a etale cyclic $m$-covering, and the composition $E \to D \to C$ is Galois with non-abelian group $\mathbb{Z}/m\mathbb{Z} \rtimes \mathbb{Z}/n\mathbb{Z}$. We will describe how to construct such non-abelian coverings. In the case of $gcd(m,n)=1$, we will also give an explicit formula counting the number of such $E$ up to isomorphism.

\section{The Main Result}

Assume the curve $C$ has genus $g \geq 2$, then by Riemann-Hurwitz formula, the genus of the curve $D$, which we denote by $h$, is $h=n(g-1)+1$. Let $\sigma$ be a generator of the Galois group of the covering $\pi: D \to C$, $\left\langle \sigma \right\rangle \cong \mathbb{Z} / n\mathbb{Z}$. Suppose the etale covering $D \to C$ corresponds to $\mathcal{L}$, which is a $n$-torsion line bundle over curve $C$. In order to search for etale cyclic $m$-covering of $D$, we should look for the $m$-torsion points in the Jacobian of $D$: $J_{D}[m]$. Over $\mathbb{C}$, $J_{D}[m] \cong (\mathbb{Z} / m \mathbb{Z})^{\oplus 2h}$. However, not all of these $m$-torsion points in $J_{D}$ would give rise to the covering $ E \to D$ we are looking for. In fact, in order to understand what the composition $ E \to D \to C$ is, we need to understand the induced action of $\sigma$ on $J_{D}[m]$.
\\ \hspace{\fill} \\
Since $\sigma$ is an isomorphism of curve $D$, the induced action of $\sigma$ on the Jacobian $J_{D}$ is also an isomorphism of $J_{D}$, which we denoted as $\sigma$ as well. In particular, it preserves $m$-torsion points in $J_{D}$, thus $\sigma $ is also an isomorphism on $J_{D}[m]$. By definition,
$J_{D} = H^{0}(D,\omega_{D})^{*} / H_{1}(D,\mathbb{Z})$. We need first understand the action of $\sigma$ on the canonical sheaf of $D$. Actually, since the covering $\pi: D \to C$ is etale, we have $\pi^{*}\omega_{C} = \omega_{D}$. Also recall that:
\begin{equation}
    \pi_{*}\mathcal{O}_{D} \cong \mathcal{O}_{C} \oplus \mathcal{L}^{-1} \oplus \cdots \oplus \mathcal{L}^{-n+1}
\end{equation}
Combine this with projection formula, we get:
\begin{equation}
    \pi_{*}\omega_{D} = \pi_{*} \pi^{*} \omega_{C} \cong \omega_{C} \otimes \pi_{*}\mathcal{O}_{D} \cong \bigoplus_{i=0}^{n-1} \omega_{C} \otimes \mathcal{L}^{-i}
\end{equation}
Taking global sections, we obtain:
\begin{equation}
    H^{0}(D,\omega_{D}) \cong \bigoplus_{i=0}^{n-1} H^{0}(C,\omega_{C} \otimes \mathcal{L}^{-i})
\end{equation}
And the right hand side of (4.3) is exactly the decomposition into $\sigma$-representations. More precisely, each $H^{0}(C,\omega_{C} \otimes \mathcal{L}^{-i})$ is an eigenspace of $\sigma$, with eigenvalue $\zeta^{i}$, where $\zeta$ is a primitive $n$-th root of unity. As one can see from (4.3), the trivial subspace of $\sigma$ inside $H^{0}(D,\omega_{D})$ is just $H^{0}(C,\omega_{C})$, and we single out the non-trivial representation:
\begin{equation}
    R:= \bigoplus_{i=1}^{n-1} H^{0}(C,\omega_{C} \otimes \mathcal{L}^{-i})
\end{equation}
By Riemann-Roch, the dimension:
\begin{equation}
    dim_{\mathbb{C}}H^{0}(C, \omega_{C})=g
\end{equation}
and for $i \neq 0$ :
\begin{equation}
    dim_{\mathbb{C}}H^{0}(C,\omega_{C} \otimes \mathcal{L}^{-i})=g-1
\end{equation}
Thus we get:
\begin{equation}
    dim_{\mathbb{C}}R=(n-1)(g-1)
\end{equation}
Here, the space $R$ is nothing but the tangent space of the Prym variety of the covering $\pi: D \to C$. By definition, the Prym variety $P$ of the covering $\pi$ is the principle connected component of the kernel of the induced norm map $Nm(\pi): J_{D} \to J_{C}$ \cite{agostini2020prym}. Alternatively, $P=Im(1-\sigma)=ker(1+\sigma + \cdots + \sigma^{n-1})^{0}$. 
\\ \hspace{\fill} \\
\textbf{Theorem 4.1:} Suppose $v \in (\mathbb{Z} / m \mathbb{Z})^{\oplus 2h} \cong J_{D}[m]$ is primitive, satisfying $\sigma v = \lambda v$ for 
some $\lambda \in (\mathbb{Z} / m \mathbb{Z})^{\times}$ and $\lambda \neq 1$. (i.e. $ v$ is an eigenvector of $\sigma: J_{D}[m] \to J_{D}[m]$, with eigenvalue $\lambda $) Then the etale cyclic $m$-covering $E \to D$ corresponding to $v$ will gives a connected non-abelian Galois covering of $C$ after composing with $\pi: D \to C$. Conversely, every connected etale Galois covering of $C$ with Galois group $\mathbb{Z}/m\mathbb{Z} \rtimes \mathbb{Z}/n\mathbb{Z}$ takes such form.
\\ \hspace{\fill} \\
\textbf{Proof:} The requirement that $v$ is primitive is to guarantee the order of $v$ is exactly $m$, which is equivalent to the requirement that the curve $E$ in covering $ E \to D$ is connected. In order to check the covering $E \to C$ is Galois, we again need to look at the action of $\sigma$ on $v$. As in our assumption, $\sigma v = \lambda v$ for some $\lambda \in (\mathbb{Z} / m \mathbb{Z})^{\times}$, and recall that $v$ and $\lambda v$ give rise to the same curve $E$, we see $\sigma$ actually lifts to an automorphism of curve $E$. Therefore, the covering $E \to C$ is Galois. To determine the Galois group, let $\theta$ be a generator of the covering $E \to D$ corresponds to eigenvector $v$, $\left\langle \theta \right\rangle \cong \mathbb{Z} / m \mathbb{Z}$. Suppose $v$ is represented by a $m$-torsion line bundle $\eta$, and
\begin{equation}
    E = Spec(\mathcal{O}_{D} \oplus \mathcal{\eta}^{-1} \oplus \cdots \oplus \mathcal{\eta}^{-m+1})
\end{equation}
Then the action of $\sigma$ sends $\eta$ to $\eta ^{\lambda}$. On each $\eta^{-i}$, the action of $\theta$ is multiply by $\epsilon^{i}$, where $\epsilon$ is a primitive $m$-th root of unity. From these, we can see that each $\eta^{k}$ is an eigen-subbundle for $\sigma \theta \sigma^{-1}$, whose action is multiplying by $\epsilon^{\lambda^{-1}k}$. Therefore, we deduce that $\sigma \theta \sigma^{-1}=\theta^{\lambda^{-1}}$. Hence, the Galois group of the covering $E \to C$ is:
\begin{equation}
    \left\langle \sigma, \theta | \sigma \theta \sigma^{-1}=\theta^{\lambda^{-1}}, \sigma^{n}=\theta^{m}=1 \right\rangle
\end{equation}
Once $\lambda \neq 1$, the group above is non-abelian, which is nothing but a semi-direct product $\mathbb{Z}/m\mathbb{Z} \rtimes \mathbb{Z}/n\mathbb{Z}$.
\\ \hspace{\fill} \\
For the converse, given a etale Galois covering $E \to C$ with Galois group $\cong \mathbb{Z}/m\mathbb{Z} \rtimes \mathbb{Z}/n\mathbb{Z}$, we consider the intermediate curve $D:= E / (\mathbb{Z}/m\mathbb{Z})$. Since $\mathbb{Z}/m\mathbb{Z}$ is a normal subgroup of its Galois group, the covering $D \to C$ is Galois and etale, with Galois group $\left\langle \sigma \right\rangle \cong \mathbb{Z}/n\mathbb{Z} $. Therefore, $D \to C$ is determined by a $n$-torsion point in $J_{C}$, and $E \to D$ is determined by a $m$-torsion point in $v \in J_{D}$. Since $E$ maps to itself under the action of $\sigma$, we deduce that $v$ must map to $\lambda v$ under the action of $\sigma$, for some $\lambda \in (\mathbb{Z} / m \mathbb{Z})^{\times}$. In order to make the covering non-abelian, it requires $\lambda \neq 1$. 
\\ \hspace{\fill} \\
Moreover, we wish to count the number of such curve $E$, given $D \to C$. Using theorem 4.1, this is equivalent to counting the number of eigen-directions (eigenvectors up to a scalar) of $\sigma$ acting on $(\mathbb{Z} / m \mathbb{Z})^{\oplus 2h}$, with eigenvalue $\neq 1$. This turns out to be a little tricky. To ease the computation, we only consider the case $gcd(m,n)=1$. We will give two methods of counting such $E$, in the section 4.1 and 4.2.

\subsection{First Method of Counting}

The idea of the first method is using similarity between matrices. To begin with, recall that the Jacobian $J_{D}=H^{0}(D,\omega_{D})^{*}/H_{1}(D,\mathbb{Z})$. Since $\sigma$ is an automorphism of $J_{D}$, it is also an isomorphism of the lattice $H_{1}(D,\mathbb{Z}) \cong \mathbb{Z}^{\oplus 2h}$. In this way, $\sigma$ becomes an element in $GL_{2h}(\mathbb{Z})$, satisfying $\sigma^{n}=id$. Since we know the eigenspace decomposition of $H^{0}(D,\omega_{D})$ under the action of $\sigma$, we can easily deduce what matrix is $\sigma$ similar to over $\mathbb{C}$. To see this, the complexification of $\sigma$, namely $\sigma \otimes \mathbb{C}$, acts on $H_{1}(D,\mathbb{Z})\otimes \mathbb{C}$, which is $H_{1}(D,\mathbb{C})$. Since $H^{1}(D,\mathbb{C}) = H^{1,0}(D)\oplus H^{0,1}(D)$, taking dual, we obtain $H_{1}(D,\mathbb{C}) \cong H^{0}(D,\omega_{D}) \oplus H^{0}(D,\omega_{D})^{*}$. Further recall that: $ H^{0}(D,\omega_{D}) \cong H^{0}(C,\omega_{C}) \oplus \bigoplus_{i=1}^{n-1} H^{0}(C,\omega_{C} \otimes \mathcal{L}^{-i})$ gives the eigenspace decomposition of the action of $\sigma$, we obtain:
\begin{equation}
    \sigma \otimes \mathbb{C} = I_{2g} \oplus \zeta I_{2g-2} \oplus \cdots 
    \oplus \zeta^{n-1} I_{2g-2}
\end{equation}
If we denote $J=diag(\zeta, \zeta^{2}, \cdots , \zeta^{n-1})$, then $\sigma$ is similar to the matrix $I_{2g} \oplus (J \otimes I_{2g-2})$ over $\mathbb{C}$.
\\ \hspace{\fill} \\
Define the $(n-1) \times (n-1)$ matrix $G=$:
$\begin{pmatrix}
    0 & \cdots & 0& -1 \\  1 & \cdots & 0 & -1\\ \vdots &\ddots & \vdots &\vdots \\ 0 &\cdots & 1 & -1 
\end{pmatrix}$
It is easy to see $G$ is similar to $J$ over $\mathbb{C}$. Thus, $\sigma$ is similar to $I_{2g} \oplus (G \otimes I_{2g-2})$ over $\mathbb{C}$. Note that both $\sigma$ and $I_{2g} \oplus (G \otimes I_{2g-2})$ over $\mathbb{C}$ is defined over $\mathbb{Z}$, thus over $\mathbb{Q}$, we deduce that $\sigma$ and $I_{2g} \oplus (G \otimes I_{2g-2})$ is similar over $\mathbb{Q}$. The difficulty here we cannot obtain similarity over $\mathbb{Z}$. To make the argument work, we must require $gcd(m,n)=1$ in the following.
\\ \hspace{\fill} \\
The goal is to obtain the similarity type of $\sigma$ after reduction modulo $m$. Once we get that information, we will know how $\sigma$ acts on $\frac{1}{m} H_{1}(D,\mathbb{Z}) /H_{1}(D,\mathbb{Z}) $, which is just the $m$-torsion points in $J_{D}$. We achieve our goal in several steps:
\\ \hspace{\fill} \\
\textbf{Step 1:} 
\\ Take any prime $p|m$, consider the similarity type of $\sigma (p) := \sigma$ mod $p$. Here we view $\sigma (p)$ as an element in $GL_{2h}(\mathbb{F}_{p})$. Because $\sigma$ is similar to $I_{2g} \oplus (G \otimes I_{2g-2})$ over $\mathbb{Q}$, they have the same characteristic polynomial and minimal polynomial over $\mathbb{Q}$. Therefore, the characteristic polynomial of $\sigma$ equals to $(x-1)^{2g}(1+x+\cdots + x^{n-1})^{2g-2}$, and the minimal polynomial of $\sigma$ equals to $x^{n}-1$. The characteristic polynomial of $\sigma (p)$ is just $(x-1)^{2g}(1+x+\cdots + x^{n-1})^{2g-2}$ mod $p$, and the minimal polynomial of $\sigma (p)$, say $q(x)$, must dividing $x^{n}-1 $ mod $p$. 
Since $gcd(p,n)=1$, $x^{n}-1$ has distinct roots over $\overline{\mathbb{F}}_{p}$. Assume $q(x) \neq x^{n}-1$ mod $p$, then some factors of $x^{n}-1$ must be left out by $q(x)$. As a result, such factors cannot appear in the characteristic polynomial $(x-1)^{2g}(1+x+\cdots + x^{n-1})^{2g-2}$ mod $p$. However, $(x-1)^{2g}(1+x+\cdots + x^{n-1})^{2g-2}$ mod $p$ has the same factor with  $x^{n}-1$ mod $p$, leading to a contradiction. Thus, $q(x)=x^{n}-1$ mod $p$. 
\\ \hspace{\fill} \\
Therefore, $\sigma(p)$ and $I_{2g} \oplus (G \otimes I_{2g-2})$ mod $p$ have the same characteristic and minimal polynomial. Note that the minimal polynomial $x^{n}-1$ does not have multiple root over $\overline{\mathbb{F}}_{p}$, they are both diagonalizable over $\overline{\mathbb{F}}_{p}$. Same characteristic polynomial gives the same multiplicity of eigenvalues. Thus, over $\overline{\mathbb{F}}_{p}$, $\sigma(p)$ is similar to  $I_{2g} \oplus (G \otimes I_{2g-2})$. Notice both matrices are defined over $\mathbb{F}_{p}$, they are actually similar over $\mathbb{F}_{p}$.
\\ \hspace{\fill} \\
\textbf{Step 2:} 
\\ In order to lift the similarity relation from over $\mathbb{F}_{p}$ to over $\mathbb{Z} / m \mathbb{Z}$, we need some results of J.Pomfret \cite{pomfret1973similarity}:
\\ \hspace{\fill} \\
\textbf{Lemma 4.2:} Let $R$ be a finite local ring with maximal ideal $M$, and $R/M=\mathbb{F}_{p^{f}}=k$. Let $\alpha$, $\beta$ be elements of $GL_{n}(R)$ with $gcd(|\left\langle\alpha \right\rangle|,p)=gcd(|\left\langle \beta \right\rangle|,p)=1$ (where $|\left\langle x \right\rangle|$ means the order of $x$). Then $\alpha$ is similar to $\beta$ if and only if $\alpha$ is similar to $\beta$ modulo $M$.
\\ \hspace{\fill} \\
If $R$ is a finite commutative ring with identity, then $R$ is uniquely isomorphic to a product of finite local rings: $R \cong \prod_{i=1}^{t} R_{i}$. Suppose $R_{i}$ is a finite local ring with maximal ideal $M_{i}$, and $R_{i}/M_{i}=k_{i}$, we have the epimorphisms:
\begin{equation}
    GL_{n}(R) \cong \prod_{i=1}^{t} GL_{n}(R_{i}) \stackrel{\pi_{i}}{\longrightarrow} GL_{n}(R_{i}) \to GL_{n}(k_{i})
\end{equation}
Combining lemma 4.2 and (4.11), we get:
\\ \hspace{\fill} \\
\textbf{Lemma 4.3:} Let $R$ be a finite commutative ring with identity and let the cardinality of $R$ be $m$. Two elements $\alpha$ and $\beta$ of $GL_{n}(R)$ satisfying $gcd(|\left\langle \alpha \right\rangle|,m)=gcd(|\left\langle \beta \right\rangle|,m)=1$ are similar if and only if they are similar over each residue field $k_{i}$.
\\ \hspace{\fill} \\
In step 1, we see $\sigma (p)$ and $I_{2g} \oplus (G \otimes I_{2g-2})$ are similar over $\mathbb{F}_{p}$ for any $p|m$. Also note that $\mathbb{F}_{p}$'s are all the residue fields of $\mathbb{Z} / m \mathbb{Z}$.
Applying lemma 4.3 to our case: $R=\mathbb{Z} / m \mathbb{Z}$, $\alpha = \sigma (m)$ ($=\sigma$ mod $m$), and $\beta = I_{2g} \oplus (G \otimes I_{2g-2}) $, we see that $\sigma (m)$ is similar to $I_{2g} \oplus (G \otimes I_{2g-2})$ over $\mathbb{Z} / m \mathbb{Z}$.
\\ \hspace{\fill} \\
\textbf{Step 3:}
\\ Denote $A = G \otimes I_{2g-2}$, viewed as an element in $GL_{(2g-2)(n-1)}(\mathbb{Z} / m \mathbb{Z})$. The part $I_{2g}$ corresponds to the eigenspace with eigenvalue $1$. Since we only want the eigenspace with eigenvalue $\lambda \neq 1$, it is equivalent to count the eigen-directions that $A$ have acting on $(\mathbb{Z} / m \mathbb{Z})^{\oplus (2g-2)(n-1)}$. 
\\ \hspace{\fill} \\
We need to understand $(\mathbb{Z} / m \mathbb{Z})^{\oplus (2g-2)(n-1)}$ as a $(\mathbb{Z} / m \mathbb{Z})[A]$ module. Since $A=G \otimes I_{2g-2}$, we only need to study the $(\mathbb{Z} / m \mathbb{Z})[G]$-mod structure on $(\mathbb{Z} / m \mathbb{Z})^{\oplus (n-1)}$, because as a $(\mathbb{Z} / m \mathbb{Z})[A]$ module, $(\mathbb{Z} / m \mathbb{Z})^{\oplus (2g-2)(n-1)}$ is just isomorphic to the direct sum of $(2g-2)$-many copies of $(\mathbb{Z} / m \mathbb{Z})^{\oplus (n-1)}$, with each summand carrying the same $(\mathbb{Z} / m \mathbb{Z})[G]$-mod structure. The minimal polynomial and the characteristic polynomial of $G$ both equals to $f(x)=1+x+x^{2}+\cdots + x^{n-1}$. Therefore, as a $(\mathbb{Z} / m \mathbb{Z})[G]$-mod, $(\mathbb{Z} / m \mathbb{Z})^{\oplus (n-1)}$ is isomorphic to $(\mathbb{Z} / m \mathbb{Z})[G]$. 
\\ \hspace{\fill} \\
Now we compute $(\mathbb{Z} / m \mathbb{Z})[G]$. We have the decomposition:
\begin{equation}
    (\mathbb{Z} / m \mathbb{Z})[G] \cong \prod_{p|m}(\mathbb{Z} / p^{e} \mathbb{Z})[G]
\end{equation}
For a set of eigenvectors: $\{v_{p} \in (\mathbb{Z} / p^{e} \mathbb{Z})^{\oplus (n-1)}\}_{p|m}$, satisfying $Gv_{p}=\lambda_{p}v_{p}$ mod $p^{e}$, we can find a unique vector $v \in (\mathbb{Z} / m \mathbb{Z})^{\oplus (n-1)} $, satifying $Gv=\lambda v$, and $v \equiv v_{p}$ mod $p^{e}$, $\lambda \equiv \lambda_{p}$ mod $p^{e}$. In this way, we reduce the computation to $(\mathbb{Z} / p^{e} \mathbb{Z})[G]$. 
\\ \hspace{\fill} \\
Note that $(\mathbb{Z} / p^{e} \mathbb{Z})[G] \cong (\mathbb{Z} / p^{e} \mathbb{Z})[x]/f(x)$, and $f(x) = \prod_{d|n, d\neq 1} \Phi_{d}(x)$. It suffices for us to consider the behavior of decomposition of each $\Phi_{d}(x)$ mod $p^{e}$. We have the following lemma:
\\ \hspace{\fill} \\
\textbf{Lemma 4.4:} For $d \nmid p-1$, $\Phi_{d}(x)$ has no degree one factor modulo $p^{e}$. For $d\mid p-1$, $\Phi_{d}(x)$ decompose completely into $\phi (d)$-many degree one factors modulo $p^{e}$.
\\ \hspace{\fill} \\
\textbf{Proof:} Passing to the residue field of $\mathbb{Z} / p^{e} \mathbb{Z}$, which is just $\mathbb{F}_{p}$, we first the decomposition of $\Phi_{d}(x)$ in $\mathbb{F}_{p}[x]$. If $\Phi_{d}(x)$ mod $p^{e}$ has any linear term $(x-\lambda)$, then so does $\Phi_{d}(x)$ mod $p$, because $\lambda \in (\mathbb{Z} / p^{e} \mathbb{Z})^{\times}$ implies $\lambda \neq 0$ mod $p$. This shows that if $\Phi_{d}(x)$ mod $p$ does not have any linear factor, then so does $\Phi_{d}(x)$ mod $p^{e}$. The only case that $\Phi_{d}(x)$ has root over $\mathbb{F}_{p}$ is that: the order of $p$ in $\mathbb{Z} / d \mathbb{Z}$ is exactly one, i.e. $d\mid p-1$. (In this case, primitive $d$-th roots of unity will be fixed under Frobenius, thus lie in $\mathbb{F}_{p}$).  
\\ Assuming $d \mid p-1$, then $\Phi_{d}(x)$ splits completely in $\mathbb{F}_{p}[x]$:
\begin{equation}
    \Phi_{d}(x) \equiv (x-a_{1})(x-a_{2}) \cdots (x-a_{\phi(d)})
\end{equation}
Since $gcd(d,p)=1$, roots $a_{i}$ and $a_{j}$ are distinct. By Hensel's lemma, we can lift the factorization into $\mathbb{Z} / p^{e} \mathbb{Z}$, namely:
\begin{equation}
    \Phi_{d}(x) \equiv (x-\Tilde{a}_{1})(x-\Tilde{a}_{2}) \cdots (x-\Tilde{a}_{\phi(d)})
\end{equation}
with $\Tilde{a}_{i} \in (\mathbb{Z} / p^{e} \mathbb{Z})^{\times}$ distinct. The lemma is proven. 
\\ \hspace{\fill} \\
With the aid of lemma 4.4, let us count how many linear factors can $f(x)$ factor out in $(\mathbb{Z} / p^{e} \mathbb{Z})[x]$. Only for those $d$, satisfying: $d \mid n$, $d \neq 1$, and $d \mid p-1$, $\Phi_{d}(x)$ will contribute $\phi (d)$-many linear factors. Therefore, the total number of linear factors in $(\mathbb{Z} / p^{e} \mathbb{Z})[x]/f(x)$ equals to:
\begin{equation}
    \sum_{d \mid gcd(n,p-1), d \neq 1} \phi (d) = gcd(n,p-1)-1
\end{equation}
Tensoring by $I_{2g-2}$, we pass from $(\mathbb{Z} / p^{e} \mathbb{Z})[G]$ to $(\mathbb{Z} / p^{e} \mathbb{Z})[A]$. By doing so, the one dimensional eigenspace with eigenvalue $\Tilde{a}_{i}$ (corresponding to the linear factor $(x-\Tilde{a}_{i})$) is enlarged into a $(2g-2)$-dimensional eigenspace still with eigenvalue $\Tilde{a}_{i}$. Since the roots of $f(x)$ are distinct, all the eigenvectors come from the union of such $(2g-2)$-dimensional eigenspaces. We first count the number of primitive points in those eigenspaces: 
\\ \hspace{\fill} \\
Take a $N$ dimensional space, $V=(\mathbb{Z} / p^{e} \mathbb{Z})^{\oplus N}$, the primitive vectors (those elements of order $p^{e}$) in $V$ can be viewed as $(\mathbb{Z} / p^{e} \mathbb{Z})^{\oplus N} - (p\mathbb{Z} / p^{e} \mathbb{Z})^{\oplus N}$. Therefore, the number of primitive vectors in $V$ equals to $p^{eN}(1-p^{-N})$. Furthermore, suppose $\mathbb{Z} / m \mathbb{Z} \cong \prod_{p|m}\mathbb{Z} / p^{e} \mathbb{Z}$, and  
suppose $X$ is a $N$ dimensional space over $\mathbb{Z} / m \mathbb{Z}$, say $X=(\mathbb{Z} / m \mathbb{Z})^{\oplus N}$, then a vector $v \in X$ is primitive if and only if $v_{p}=v$ mod $p^{e}$ is primitive in $X_{p} = (\mathbb{Z} / p^{e} \mathbb{Z})^{\oplus N}$ for every $p \mid m$. 
\\ \hspace{\fill} \\
Combining all the arguments above together, we obtain the number of primitive eigenvector for $A$ in $(\mathbb{Z} / m \mathbb{Z})^{\oplus (2g-2)(n-1)}$ equals to:
\begin{equation}
    \prod_{p \mid m} (gcd(n,p-1)-1)\cdot p^{e(2g-2)}(1-p^{-2g+2})
\end{equation}
Since what we need to count is the number of eigen-directions, we need to modulo thus primitive eigenvectors by scalars in $(\mathbb{Z} / m \mathbb{Z})^{\times}$. Thus, we need to divide (4.16) by $\phi (m)$:
\begin{equation}
    T_{g}(m,n)=\prod_{p \mid m} (gcd(n,p-1)-1)\cdot p^{e(2g-2)}(1-p^{-2g+2})/p^{e}(1-p^{-1})
\end{equation}
Where $T_{g}(m,n)$ is just the number we are looking for. Reduce (4.17), we obtain: 
\begin{equation}
    T_{g}(m,n)=m^{2g-3} \prod_{p \mid m} (gcd(n,p-1)-1)\cdot (1+p^{-1}+ \cdots + p^{-2g+3})
\end{equation}
Summing up, we get the following counting result:
\\ \hspace{\fill} \\
\textbf{Theorem 4.5:} Let $C$ be a smooth projective curve of genus $g \geq 2$, and suppose $m$ and $n$ are coprime integers. For a given etale cyclic $n$-covering $\pi: D
\to C$, up to isomorphism, there are $T_{g}(m,n)$-many curves $E$, such that $E \to D$ is an etale cyclic $m$-covering, and the composition $E \to C$ is Galois with non-abelian Galois group $\mathbb{Z}/m\mathbb{Z} \rtimes \mathbb{Z}/n\mathbb{Z}$.

\subsection{Second Method of Counting}

Here we will give an alternative method to count the number of such curves $E$. Recall that the $m$-torsion points in $J_{D}[m]$ are nothing but $\frac{1}{m} H_{1}(D,\mathbb{Z}) /H_{1}(D,\mathbb{Z}) $. If we can spell out the action of $\sigma$ on $H_{1}(D,\mathbb{Z})$, we will get an explicit matrix representation of $\sigma \in GL_{2h}(\mathbb{Z})$, up to conjugation over $\mathbb{Z}$. That will be enough for us to compute the eigenvectors of $\sigma(m)=\sigma$ mod $m$. The key point is a topological lemma below:
\\ \hspace{\fill} \\
\textbf{Lemma 4.6:} Suppose $\Sigma_{g}$ is a closed topological surface of genus $g \geq 2$, and $\Sigma_{h}$, where $h=n(g-1)+1$, is a normal covering surface of $\Sigma_{g}$: $\Sigma_{h} \to \Sigma_{g}$, with group of deck transformation $\mathbb{Z} / n \mathbb{Z}$. Then such covering is equivalent to the 'standard' cyclic covering between surfaces, as illustrated by Figure 1 \cite{merkulov2003hatcher} and Figure 2 \cite{lanier2018normal}.
\begin{figure}
\centering
\includegraphics[width=0.35\linewidth]{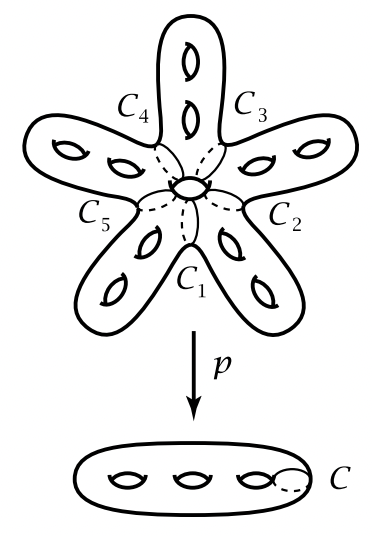}
\caption{\label{fig: cyclic1} Standard cyclic covering for $\mathbb{Z} / 5 \mathbb{Z}$}
\end{figure}
\begin{figure}
\centering
\includegraphics[width=0.6\linewidth]{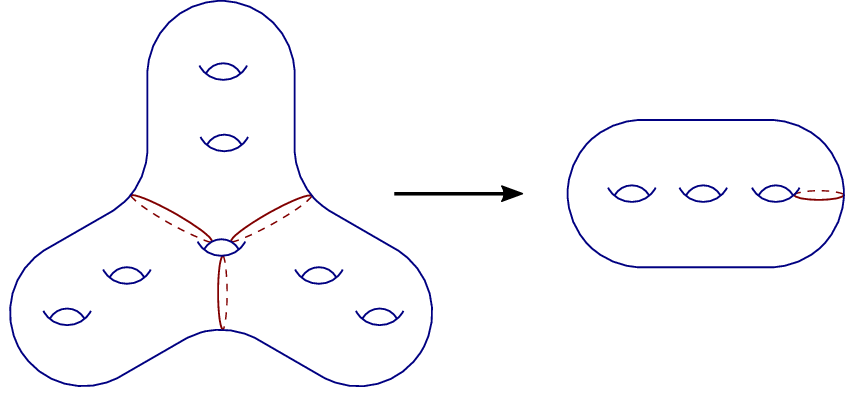}
\caption{\label{fig: cyclic2} Standard cyclic covering for $\mathbb{Z} / 3 \mathbb{Z}$}
\end{figure}
\\ \hspace{\fill} \\
We first need to translate the word standard into more concrete definition. Let $\pi = \pi_{1}(\Sigma_{g})$ be the fundamental group of genus $g$ closed surface, and let $a_{1}, b_{1}, \cdots, a_{g}, b_{g}$ be a basis of $\pi$ satisfying $\prod_{i=1}^{g}[a_{i},b_{i}]=e$. Then the standard $\mathbb{Z} / n \mathbb{Z}$-covering corresponds to the epimorphism $\pi \to \mathbb{Z} / n \mathbb{Z} \to 1$, sending $a_{1}$ to a generator of $\mathbb{Z} / n \mathbb{Z}$ and $b_{1}, a_{2}, b_{2}, \cdots a_{g}, b_{g}$ to $0$. In this case, we have a short exact sequence:
\begin{equation}
    1 \to K \to \pi \to \mathbb{Z} / n \mathbb{Z} \to 1
\end{equation}
and $K$ is nothing but the fundamental group of $\Sigma_{h}$, with $h=n(g-1)+1$. It is easy to see that such covering $\Sigma_{h} \to \Sigma_{g}$ has the geometric figure that looks like rotating by $\mathbb{Z} / n \mathbb{Z}$ along a hole, which we refer to as standard cyclic covering. Therefore, to prove Lemma 4.6, it suffices to prove the following result:
\\ \hspace{\fill} \\
\textbf{Lemma 4.7:} Let $\pi$ be the fundamental group of $\Sigma_{g}$.
For every epimorphism $p: \pi \to \mathbb{Z} / n \mathbb{Z}$, there is a basis of $\pi$, namely $\Tilde{x}_{1}, \cdots, \Tilde{x}_{g}, \Tilde{y}_{1}, \cdots, \Tilde{y}_{g}$, such that under $p$, $\Tilde{x}_{1}$ maps to generator of $\mathbb{Z} / n \mathbb{Z}$, and $\Tilde{x}_{2}, \cdots, \Tilde{x}_{g}, \Tilde{y}_{1}, \cdots, \Tilde{y}_{g}$ map to $0$. Moreover, $\prod_{i=1}^{g}[\Tilde{x}_{i},\Tilde{y}_{i}]=e$.
\\ \hspace{\fill} \\
\textbf{Proof:} Consider the following diagram:
\begin{displaymath}
    \xymatrix{ 
        \quad & \mathbb{Z} \ar[d]^{mod \quad n} & \quad \\
        \pi / [\pi, \pi] \cong \mathbb{Z}^{2g} \ar[r]_{\quad \delta} \ar[ur]^{\Tilde{\delta}} & \mathbb{Z} / n \mathbb{Z} \ar[r] \ar[d] & 0 \\
        \quad & 0 & \quad
    }
\end{displaymath}

Since $\mathbb{Z} / n \mathbb{Z}$ is abelian, $p: \pi \to \mathbb{Z} / n \mathbb{Z}$ factors through $\pi / [\pi,\pi] = H_{1}(\Sigma_{g},\mathbb{Z})$, which is isomorphic to $\mathbb{Z}^{2g}$. Since $\mathbb{Z}^{2g}$ is a free $\mathbb{Z}$-mod, thus projective, $\delta$ lifts to $\Tilde{\delta}: \mathbb{Z}^{2g} \to \mathbb{Z}$. There exists $v \in \mathbb{Z}^{2g}$, such that $\Tilde{\delta}(v)$ mod $n$ = $\delta (v)$ generates $\mathbb{Z} / n \mathbb{Z}$. If $v$ is not primitive, we can divide $v$ by $k \in \mathbb{Z}$, such that $v/k$ is primitive and $\Tilde{\delta}(v) \in \mathbb{Z}$. After modulo $n$, $\delta(v/k)$ still generates $\mathbb{Z} / n \mathbb{Z}$. Therefore, we can find a primitive vector $v \in \mathbb{Z}^{2g}$, satisfying $\delta (v) \in (\mathbb{Z} / n \mathbb{Z})^{\times}$. Now since $v$ is primitive, $\mathbb{Z}^{2g} = \mathbb{Z}\{ v \} \oplus ker(\Tilde{\delta})$, and $ker(\Tilde{\delta}) \cong \mathbb{Z}^{2g-1}$.
\\ \hspace{\fill} \\
Consider the intersection form on $H_{1}(\Sigma_{g},\mathbb{Z})$, i.e. $B$, a skew-symmetric, bilinear, unimodular form on $\mathbb{Z}^{2g}$. Note that $B(v,v)=0$. Consider the linear function $B(v,\cdot): ker(\Tilde{\delta}) \to \mathbb{Z}$. It must be surjective, otherwise $det(B) \neq \pm 1$. Therefore, exists $y_{1} \in ker(\Tilde{\delta})$, primitive, such that $\{ v, y_{1} \}$ form a hyperbolic basis. $ ker(\Tilde{\delta)} \cong \mathbb{Z}^{2g-1}  = \mathbb{Z} \{ y_{1} \} \oplus U$, with $\mathbb{Z}^{2g-2} \cong U=ker(\Tilde{\delta}) $.
\\ \hspace{\fill} \\
Further consider the linear function $B(y_{1},\cdot): U \to \mathbb{Z}$. If $ker(B(y_{1},\cdot))$ is the whole $U$, then $B$ takes diagonal form under $ \mathbb{Z}\{ v, y_{1} \} \oplus U$, and the restriction of $B$ to $U$ is also unimodular and skew-symmetric, and we can apply induction to get a standard basis. If $B(y_{1},\cdot)$ is non-trivial, let $W \subset U$ to be the kernel of $B(y_{1},\cdot)$. Note that $W \cong \mathbb{Z}^{2g-3}$, and $B$ restricts to $W$ is integral and skew-symmetric. Thus, there is a basis $\{ w_{1}, w_{2}, \cdots, w_{2g-3} \}$ for $W$. We can assume that $B$ takes the form $diag(S_{1}, S_{2}, \cdots, S_{g-2}, 0)$ on $W$, where $S_{i}$'s are integral skew-symmetric $2 \times 2$ matrices, say  $S_{i}=\begin{pmatrix}
    0 & d_{i} \\
    -d_{i} & 0
\end{pmatrix}$ 
Assume $U = \mathbb{Z} \{ t \} \oplus W$, then under the basis $\{ v, y_{1}, w_{1}, \cdots, w_{2g-3}, t \} $, $B$ takes the form: 
\begin{equation}
    \begin{pmatrix}
    0 & 1 & 0 & 0 & \cdots & 0 & 0 \\
    -1 & 0 & 0 & 0 & \cdots & 0 & p_{2g-2} \\
    0 & 0 & S_{1} & 0 & \cdots & 0 & p_{2g-3} \\
    0 & 0 & 0 & S_{2} & \cdots & 0 & p_{2g-4} \\
    \vdots & \vdots & \vdots & \vdots & \ddots & \vdots & \vdots \\
    0 & 0 & 0 & 0 & \cdots & S_{g-2} & p_{2} \\
    0 & 0 & 0 & 0 & \cdots & 0 & p_{1} \\
    0 & -p_{2g-2} & -p_{2g-3} & -p_{2g-4} & \cdots & -p_{1} & 0
\end{pmatrix}
\end{equation}
Computing $det(B)$ from (4.20), we get:
\begin{equation}
    det(B) = d_{1}^{2}d_{2}^{2}\cdots d_{g-2}^{2}p_{1}^{2} \nonumber
\end{equation}
But $B$ is unimodular, so $d_{1}^{2}d_{2}^{2}\cdots d_{g-2}^{2}p_{1}^{2} = 1$, which implies $d_{i} = \pm 1$ and $p_{1}= \pm 1$. As a result, we are able to perform changing basis within sublattice $\mathbb{Z}^{2g-1} = \mathbb{Z} \{ y_{1}, w_{1},\cdots, w_{2g-3}, t \} = ker(\Tilde{\delta})$, to make $B$ into the standard form:
\begin{equation}
    B \sim diag\{ \begin{pmatrix}
        0 & 1 \\
        -1 & 0
    \end{pmatrix}, \cdots, \begin{pmatrix}
        0 & 1 \\
        -1 & 0
    \end{pmatrix} \} \nonumber
\end{equation}
Summing up, we can find a symplectic basis $\{ x_{1}, \cdots, x_{g}, y_{1}, \cdots , y_{g} \}$ of $\mathbb{Z}^{2g}$, such that the intersection form $B$ is represented by:
\begin{equation}
    B= x_{1}^{*} \wedge y_{1}^{*} + \cdots + x_{g}^{*}\wedge y_{g}^{*} \nonumber
\end{equation}
Moreover, $\delta (x_{1})$ generates $\mathbb{Z} / n \mathbb{Z}$, and $\delta(x_{2}), \cdots, \delta(x_{g}), \delta(y_{1}), \cdots, \delta(y_{g}) =0$. By realization theorem, we can also find $\Tilde{x_{i}}, \Tilde{y_{i}} \in \pi = \pi_{1}(\Sigma_{g})$, such that $\Tilde{x_{i}} \to x_{i}$, and $\Tilde{y_{i}} \to y_{i}$, under the Hurewicz map $\pi_{1}(\Sigma_{g}) \to H_{1}(\Sigma_{g}, \mathbb{Z})$, and the geometric intersection number 
$ \#(\Tilde{x_{i}}, \Tilde{y_{j}})=\delta_{ij} $. 
\\ \hspace{\fill} \\
Let $\{ a_{i}, b_{j} \}$ be the standard basis in $H_{1}(\Sigma_{g}, \mathbb{Z})$, such that exist $\{ \alpha_{i}, \beta_{j} \} \subset \pi_{1}(\Sigma_{g}) \to \{ a_{i}, b_{j} \} $ with relation $\prod_{i=1}^{g}[\alpha_{i}, \beta_{j}] = e $. Then the basis $\{ x_{i}, y_{j} \}$ can be obtained from $\{ a_{i}, b_{j} \}$ by a symplectic transformation. Since the mapping class group $Mod_{g}$ maps surjectively to the symplective group $Sp_{2g}(\mathbb{Z})$ \cite{putmansymplectic}, \cite{lanier2018normal}, we can find $f$, a homeomorphism of $\Sigma_{g}$, such that $f_{*}(a_{i})=x_{i}$, and $f_{*}(b_{j})=y_{j}$ (here $f_{*}: H_{1}(\Sigma_{g},\mathbb{Z}) \to H_{1}(\Sigma_{g},\mathbb{Z})$). Let $\Tilde{x_{i}}=f_{*}(\alpha_{i})$, and $\Tilde{y_{j}}=f_{*}(\beta_{j})$ (here $f_{*}: \pi \to \pi$, and $\Tilde{x_{i}}, \Tilde{y_{j}} \in \pi$). Then we have $\prod_{i=1}^{g}[\Tilde{x_{i}},\Tilde{y_{i}}]=e$, and $\Tilde{x_{i}} \to x_{i}, \Tilde{y_{j}} \to y_{j}$ in $H_{1}(\Sigma_{g},\mathbb{Z})$. So these $\{ \Tilde{x_{i}}, \Tilde{y_{j}} \}$ is just the basis we wanted, proving the lemma.
\\ \hspace{\fill} \\
Under the standard cyclic covering, the matrix for the deck transformation is $I_{2} \oplus \begin{pmatrix}
    \quad & 1 \\
    I_{n-1} & \quad
\end{pmatrix} \otimes I_{2g-2}$, acting on $H_{1}(\Sigma_{h},\mathbb{Z})$. Applying this to our case $\pi: D \to C$, we obtain:
\\ \hspace{\fill} \\
\textbf{Corollary 4.7:} $\sigma$ is similar to $I_{2} \oplus \begin{pmatrix}
    \quad & 1 \\
    I_{n-1} & \quad
\end{pmatrix} \otimes I_{2g-2}$ over $\mathbb{Z}$.
\\ \hspace{\fill} \\
Since we obtain similarity type over $\mathbb{Z}$, we can go further than we do in the first method. In fact, for any integer $m$, not necessarily coprime with $n$, by Corollary 4.7, we have $\sigma (m)$ is similar to $I_{2} \oplus \begin{pmatrix}
    \quad & 1 \\
    I_{n-1} & \quad
\end{pmatrix} \otimes I_{2g-2}$
over $\mathbb{Z} / m \mathbb{Z}$. From here on, the computation is almost the same as the third step of the first method. However, if $gcd(m,n) > 1$, the case will get more complicated, and we won't do the computation for it here.

\section{Some Other Results}
In this last section, we use the main result in section 4 to deduce some further corollaries. 
\\ \hspace{\fill} \\
We keep the assumption $gcd(m,n)=1$ throughout this section. To begin with, note that in the semi-direct product $\mathbb{Z}/m\mathbb{Z} \rtimes \mathbb{Z}/n\mathbb{Z}$, the subgroup $\mathbb{Z} / m \mathbb{Z}$ is unique. Therefore, for the covering $E \to C$, there is a unique intermidiate curve $D := E / (\mathbb{Z} / m \mathbb{Z})$, such that $D \to C$ is an etale cyclic $n$-covering. In other words, if we start our construction from non-isomorphic $D$, ($D \to C$ is etale $n$-cyclic), we will obtain non-isomorphic curve $E$. 
\\ \hspace{\fill} \\
It is easy for us to count: given a smooth projective curve $C$, the number of curve $D$ up to isomorphism, such that $D \to C$ is etale cyclic $n$-covering. From section 2.2, we know it equals the number of primitive points in $(\mathbb{Z} / n \mathbb{Z})^{\oplus 2g}$, modulo scalars in $(\mathbb{Z} / n \mathbb{Z})^{\times}$. Again, assume that $\mathbb{Z} / n \mathbb{Z} = \prod_{q \mid n} \mathbb{Z} / q^{f} \mathbb{Z}$. As in section 4.1, the number we are looking for is:
\begin{equation}
    \prod_{q \mid n} q^{2gf}(1-q^{-2g})/ q^{f}(1-q^{-1})
\end{equation}
(5.1) equals to:
\begin{equation}
    n^{2g-1} \prod_{q \mid n} (1+q^{-1}+ \cdots + q^{-2g+1})
\end{equation}
Combining these together, we obtain:
\\ \hspace{\fill} \\
\textbf{Corollary 5.1:} For a given smooth projective curve $C$ of genus $g \geq 2$, the number of curve $E$ up to isomorphism, such that $E \to C$ is etale non-abelian Galois with Galois group a semi-direct product $\mathbb{Z}/m\mathbb{Z} \rtimes \mathbb{Z}/n\mathbb{Z}$ equals to $C_{g}(m,n)$:
\begin{equation}
    C_{g}(m,n)=m^{2g-3}n^{2g-1} \prod_{p \mid m , q \mid n}(gcd(n,p-1)-1)(1+p^{-1}+ \cdots + p^{-2g+3})(1+q^{-1}+ \cdots + q^{-2g+1})
    \nonumber
\end{equation}
\\ \hspace{\fill} \\
The most special case of Corollary 5.1 might be $m=p$, $n=q$, with $p \neq q$ are distinct prime numbers. In the case of $q \nmid p-1$, $C_{g}(p,q)=0$; If $q \mid p-1$, then $C_{g}(p,q)=(q^{2g}-1)(1+p+ \cdots + p^{2g-3})$. We get:
\\ \hspace{\fill} \\
\textbf{Corollary 5.2:} Given curve $C$ as above. If $q \mid p-1$, there are $(q^{2g}-1)(1+p+ \cdots + p^{2g-3})$ many curve $E$ up to isomorphism such that $E \to C$ is etale Galois with Galois group $\mathbb{Z}/p\mathbb{Z} \rtimes \mathbb{Z}/q\mathbb{Z}$.

\bibliographystyle{alpha}
\bibliography{sample}

\end{document}